\documentclass[12pt]{article}
\usepackage{fullpage,latexsym,amsfonts,amsthm,amsmath,amscd,amssymb,color,appendix}
\usepackage[dvips]{graphicx}
\usepackage[all]{xy}
\usepackage{comment}

\usepackage{epic}

\newtheorem{lemma}{Lemma}[section]
\newtheorem{thm}[lemma]{Theorem}

\newtheorem{prop}[lemma]{Proposition}
\newtheorem{cor}[lemma]{Corollary}

\newtheorem{defn}[lemma]{Definition}

\newcommand{\ad}{\operatorname{ad}}
\newcommand{\Ad}{\operatorname{Ad}}

\def\C{{\mathbb C}}
\def\R{{\mathbb R}}

\def\N{{\mathbb N}}

\def\G{{\cal G}}

\def\sl{{\mathfrak{sl}}}
\def\su{{\mathfrak{su}}}

\def\u{{\mathfrak{u}}}
\def\k{{\mathfrak k}}

\def\z{{\mathfrak z}}
\def\a{{\mathfrak a}}
\def\h{{\mathfrak h}}

\def\so{{\mathfrak {so}}}
\def\s{{\mathfrak s}}
\def\sp{{\mathfrak{sp}}}
\def\e{{\mathfrak e}}

\def\Symm{{S}}

\def\g{\mathfrak{g}}

\def\h{{\mathfrak h}}

\def\t{\tilde}
\def\t{\mathfrak{t}}

\def\rk{{\rm rk}}

\def\R{{\mathbb R}}

\newcommand{\bea}{\begin{eqnarray}}
\newcommand{\eea}{\end{eqnarray}}

\newcommand{\Tr}{\textrm{Tr}}

\def\pr{{\rm pr}}
\def\sign{{\rm sign}}
\def\id{{\rm id}}

\DeclareMathAlphabet{\mathpzc}{OT1}{pzc}{m}{it}

\begin{document}
\title{
\begin{flushright} \small
UUITP-54/21
 \end{flushright}
\bigskip
\bigskip
Nijenhuis tensor and invariant polynomials}

\author{F. Bonechi\footnote{\small INFN Sezione di Firenze, email: francesco.bonechi@fi.infn.it} ,
J. Qiu \footnote{\small Matematiska institutionen, Institutionen för fysik och astronomi, Uppsala Universitet email: jian.qiu@math.uu.se},
M. Tarlini\footnote{\small INFN Sezione di Firenze,  email: marco.tarlini@fi.infn.it},
E. Viviani\footnote{\small INFN Sezione di Firenze and Dipartimento di Fisica, Universit\`a di Firenze, email: emanuele.viviani@unifi.it}}

\date{}

\maketitle

\begin{abstract}
We discuss the diagonalization problem of the Nijenhuis tensor in a class of {\it Poisson-Nijenhuis} structures defined on compact hermitian symmetric spaces. We study its action on the ring of invariant polynomials of a Thimm chain of  subalgebras. The existence of $\phi\,$-{\it minimal representations} defines a suitable basis of invariant polynomials that completely solves the diagonalization problem. We prove that such representations exist in the classical cases AIII, BDI, DIII and CI, and do not exist
in the exceptional cases EIII and EVII. We discuss a second general construction that in these two cases computes partially the spectrum and hints at a different behavior with respect to the classical cases.
\end{abstract}

\thispagestyle{empty}

\section{Introduction}
The notion of symplectic groupoid was introduced by A. Weinstein (\cite{Weinstein}) with the problem of quantization of the underlying Poisson manifold in mind. The basic idea is that a proper quantization must be compatible with the additional groupoid structure so that the output of the procedure is an algebra, regarded as the algebra of operators. If the quantization scheme is given by geometric quantization, both prequantization and polarization should be compatible with the groupoid structures. In \cite{WeXu} it has been shown that the natural notion of compatible prequantization  is encoded in a central extension of the symplectic groupoid; moreover, if the groupoid is prequantizable as a symplectic manifold, then such compatible prequantization always exists and is unique.
In \cite{Hawkins}, where this approach has been revived, a natural notion of {\it multiplicative polarization} has been introduced. Finding polarizations that make geometric quantization work is in general highly demanding, so that there are basically only two big classes of symplectic manifolds where this program can be completed: cotangent and Kahler manifolds. The corresponding polarizations are not in general compatible with the groupoid structure so that this ambitious picture is bit meagre in terms of concrete examples.

In \cite{Bonechi} it was shown that a non degenerate {\it symplectic Poisson-Nijenhuis structures} (PN), a particular example of bihamiltonian geometry, is a source of such multiplicative polarizations. A symplectic $PN$ structure on a smooth manifold $M$ consists of a symplectic structure $\omega$ and a Poisson structure $\pi$ such that $\omega^{-1}$ and $\pi$ are compatible ({\it i.e.} $[\pi,\omega^{-1}]=0$). Among the consequences of this property, the Nijenhuis tensor $N=\pi\circ \omega$ has vanishing torsion and defines a hierarchy of compatible Poisson structures $N^n\circ\omega^{-1}, n=0,1,\ldots$. When $N$ has maximal rank, the hamiltonian forms
\begin{equation}\label{hamiltonian_forms}
\Omega^1_{ham}=\{\alpha\in\Omega^1(M), d\alpha=d_N\alpha=0\}\;,
\end{equation}
where $d_N$ is the algebroid differential,
define a lagrangian polarization that is in general very singular and so unfit for standard geometric quantization. As explained in \cite{LGSX}, the PN structure can be integrated to a {\it multiplicative PN structure} $N_\G$ on the symplectic groupoid $\G$ integrating any Poisson structure of the hierarchy so that the polarization defined by $N_\G$ is multiplicative too. The polarization defined by the hamiltonian forms is still singular but it allows one to define the topological groupoid of Lagrangian leaves and, upon regularity assumptions, the groupoid of {\it Bohr-Sommerfeld leaves}. One can then consider its convolution algebra as the quantization. This procedure was worked out in \cite{BCQT} for $\C P_n$ where the {\it Bruhat-Poisson structure} (\cite{LuWe}) $\pi$ is compatible with the Fubiny-Study symplectic form $\omega$. The $C^*$-algebra of
the groupoid of Bohr-Sommerfeld leaves coincides with the $C^*$-algebra of quantum homogeneous spaces, as shown in \cite{Sheu}. This is a particular case of a class of $PN$ structures defined on compact hermitian symmetric spaces $M_\phi$, introduced in \cite{KRR},
where $\pi$ is the Bruhat-Poisson structure and $\omega$ the KKS symplectic form.

Motivated by the problem of quantization, the study of these PN geometries was started in \cite{BQT}, where the diagonalization of the Nijenhuis tensor was solved for the classical cases (AIII, BI, DIII, CI).
Let $M_\phi=K/K_\phi$ be a compact hermitian symmetric space, where $K$ is a compact simple Lie group, $\k$ its Lie algebra and $K_\phi\subset K$ the subgroup integrating $\k_\phi\subset\k$ the Lie subalgebra defined by the non compact root $\phi$. In particular it was proved that $i$) if $\lambda$ is an eigenvalue of any solution ${\cal M}$ of the matrix equation, called the {\it master equation},
\begin{equation}\label{master_equation}
Nd {\cal M} = d{\cal M}^-{\cal M} + {\cal M} d{\cal M}^++r d{\cal M}\,,\;\;\;d{\cal M}^++ d{\cal M}^- = kd{\cal M}\,
\end{equation}
for some $k,r\in\C$ then $\tilde{\lambda}=k\lambda +r$ is an eigenvalue of $N$; $ii$) the $\k$-moment map evaluated in a $\phi\,$-minimal representation (Definition \ref{phi_minimal_rep}) satisfies the master equation. Such representations were provided for the classical cases. By a case by case analysis, it was proved that for each of the classical cases there exists a chain of subalgebras 
$$
 \k\supset \k_1\supset \k_2\ldots\supset \k_n=\t
$$
where $\t$ is a Cartan subalgebra, such that the corresponding minors of the moment map in the $\phi\,$-minimal representation satisfy (\ref{master_equation}). We refer to it as Thimm chain of subalgebras (\cite{GS1}) The exceptional $EIII, EVII$ spaces were not considered.

In this paper we discuss a new approach to this problem.
We focus more on the structure of the Hamiltonian forms $\Omega^1_{ham}(M_\phi)$ rather than the eigenvalues of $N$ and describe them in terms of invariant polynomials with respect to the Thimm chain. Notice that Nijenhuis eigenvalues are in general only continuous functions and their derivatives have singularities. Since the multiplicative polarization of the groupoid integrates the hamiltonian forms on $M_\phi$, it is important that we describe them in terms of global variables. Moreover, we start the study of the exceptional cases that were not discussed in \cite{BQT}.

Let us describe in some details our results. We consider the subcomplex $(\Omega^{\k_1}_{pol},d)$ of the de Rham complex generated by the invariant polynomials $S(\k_1)^{\k_1}$ of any subalgebra $\k_1\subset\k$ satisfying  the compatibility condition (\ref{hypothesys_subalgebra}) with the complex structure $J$ of $\k$. We compute in Proposition \ref{thm_basic_forms} the general formula (\ref{action_Nijenhuis}) for $d_Np$ where $p\in S(\k_1)^{\k_1}$ and we write the sufficient condition (\ref{sufficient_condition}) that implies that $d_Np$ is $\k_1$-basic.
This condition is satisfied if $\k_1=\k_\phi$ for all compact hermitian symmetric spaces so that $d_Np$ is a $\k_\phi$-basic form (Proposition \ref{general_corollary}).

We consider Thimm chains of subalgebras satisfying (\ref{hypothesys_subalgebra}) in Section \ref{thimm_section}.
We are able to prove that the sufficient condition (\ref{sufficient_condition}) is satisfied for each subalgebra provided a $\phi\,$-minimal representation of $\k$ exists and we compute explicitly the action of $d_N$ on a suitable basis of invariant polynomials (Theorem \ref{explicit_formula}). In particular, in this case $(\Omega^{\k_i}_{pol},d_N)$ is a subcomplex of the Nijenhuis complex and the diagonalization is easily obtained (Corollary \ref{cor_indecomposable}). These $\phi\,$-minimal representations exist for the classical cases: in particular the diagonalization results of \cite{BQT} are better understood in this more conceptual Lie theoretical framework.

We prove also that these representations do not exist for EIII and EVII (Proposition \ref{nogo_min_rep_exceptional}) so that
Corollary \ref{cor_indecomposable} does not apply to these cases. We can still use the general result of Proposition \ref{general_corollary} that implies that $d_NS(\k_\phi)^{\k_\phi}$ are $\k_\phi$-basic forms, where $\k_\phi=\so(10)\oplus \so(2)$ for EIII and $\k_\phi=\e_6\oplus \so(2)$ for EVII. We give its explicit description in Section \ref{exceptional} and prove that $(\Omega^{\k_\phi}_{pol},d_N)$ is a subcomplex of the Nijenhuis complex for EIII. In the case of EVII a new phenomenon occurs that suggests that we have to consider a polynomial ring of invariants bigger than  $S(\k_\phi)^{\k_\phi}$ . These results show that the exceptional cases behave quite differently with respect to the classical cases and open the way to the complete solution that we plan to address in a future paper.

\section{Compact Hermitian symmetric spaces}\label{conventions}

We fix here notations and basic facts of compact hermitian symmetric spaces (see \cite{Wolf,Knapp}.
Let $\g$ be a complex semisimple Lie algebra and let $\k$ be its compact real form. We denote with $(,)$ the Killing form. Let us fix a compact Cartan subalgebra $\t\subset\k$ and let us denote by $\t_\C$ its complexification. We denote with $\Delta$,
$\Delta^+$ and $\Pi$ the roots, a choice of positive roots and the simple roots, respectively. For each $\alpha\in\Delta$ we
denote with $\g_\alpha$ the root space and we fix a root vector $e_\alpha$. We define $J:\g\rightarrow\g$ as
\begin{equation}\label{complex_structure}
J(\t)=0\,,\;\;\;\; J(e_\alpha) = i\  {\sign}(\alpha)e_\alpha\;.
\end{equation}
We fix the normalization of the root vectors as $(e_\alpha,e_{-\alpha})(\alpha,\alpha)=2$ so that $[e_\alpha,e_{-\alpha}] = h_\alpha$, where $h_\alpha$ is the coroot vector of $\alpha$.

A simple root $\phi$ is called non compact (see \cite{Wolf}) if the decomposition
of any $\alpha\in\Delta^+$ along $\phi$ is either $0$ or $1$.
Roots (and positive roots) are decomposed accordingly as $\Delta=\Delta_c\coprod\Delta_{nc}$ where $\Delta_c$ are the compact roots and $\Delta_{nc}$ the non compact ones. Since the sum of two  compact roots, if it is a root, is compact, then
\begin{equation}\label{main_subalgebra}\k_\phi = \t \dot{+}_{\alpha\in\Delta_{c}^+}(\g_\alpha\dot{+}\g_{-\alpha})\cap\k
\end{equation}
is a subalgebra. Its center $\z_\phi$ is one dimensional and is generated by $\rho_\phi\in\t$ defined as $\alpha(\rho_\phi)=0$ for each $\alpha\in\Delta_c$
and normalized as $\phi(\rho_\phi)=i$. We get the decomposition
\begin{equation}\label{decomposition_subalgebra}
\k_\phi = \k_\phi'\oplus \z_\phi
\end{equation}
with $\k_\phi'$ semisimple. We denote with $\k_\phi^\perp=\dot{+}_{\alpha\in\Delta_{nc}^+}(\g_\alpha\dot{+}\g_{-\alpha})\cap\k$ the orthogonal complement with respect to the Killing form. We have that
\begin{equation}\label{complex_structure_orthogonal}J|_{\k_\phi^\perp}= [\rho_\phi,-]\,,\;\;\;(J|_{\k_\phi^\perp})^2=-1\;.\end{equation}

If $K$ and $K_\phi\subset K$ integrate $\k$ and $\k_\phi$, we denote with $M_\phi=K/K_\phi$ the corresponding homogeneous space, that is a compact hermitian symmetric space. After identifying $\k$ with $\k^*$, $M_\phi$ can be realized as an adjoint orbit $M_\phi = K\rho_\phi K^{-1}$. We denote with $\omega$ the Konstant-Kirillov-Souriau symplectic form and with $\mu:M_\phi\rightarrow\k$ the moment map for the hamiltonian $K$-action. We will denote with $X^\sharp$ the fundamental vector field
associated to $X\in\k$ and with $v_f=\omega^{-1}df$ the hamiltonian vector field of $f\in C^\infty(M_\phi)$.

Let us consider the non compact real form
$$\g_0=\k_\phi\dot{+}i\k_\phi^\perp$$
of $\g$ where the Cartan involution $\theta$ is defined as $\theta|_{\k_\phi}=\id$ and $\theta_{\k_\phi^\perp}=-\id$. A $\theta$-stable Cartan subalgebra $\h$ of $\g_0$ is decomposed as $\h=\t_0\oplus i\a$; it is maximally non compact when the dimension of $\a$ is as large as possible. A maximal abelian subalgebra $\a\subset \k_\phi^\perp$ is unique up to $K_\phi$-conjugation (see \cite{Knapp}). We refer to $\dim \a$ as the rank of $M_\phi$.

\medskip
Here we list the classification of compact hermitian symmetric spaces.

\begin{itemize}
\item[$AIII$)] $\k=\su(n+1)$ and $\phi = \alpha_i$, $\k_\phi=\s(\u(i)\oplus \u(n+1-i))$. The rank is ${\rm min}\{i,n+1-i\}$.
\item[$BDI$)] $\k=\so(n+2)$ and $\phi = \alpha_1$, $\k_\phi=\so(n)\oplus\so(2)$. The rank is $2$.
\item[$DIII$)] $\k=\so(2n)$ and $\phi = \alpha_n$, $\k_\phi=\u(n)$. The rank is $[n/2]$.
\item[$CI$)] $\k=\sp(n)$, $\phi=\alpha_n$ and $\k_\phi=\u(n)$. The rank is $n$.
\item[$EIII$)] $\k=\e_6$, $\phi=\alpha_6$ and $\k_\phi=\so(10)\oplus\so(2)$. The rank is $2$.
\item[$EVII$)] $\k=\e_7$, $\phi=\alpha_1$ and $\k_\phi=\e_6\oplus \so(2)$. The rank is $3$.
\end{itemize}

\begin{figure}[ht]
\begin{picture}(10,50)(-20,-30)
\put(0,0){$A_n$}
\multiput(35,0)(20,0){5}{\circle{8}}
\multiputlist(45,0)(20,0)
{{\line(1,0){12}},{\dashline{2}(-5,0)(5,0)},{\dashline{2}(-5,0)(5,0)},{\line(1,0){12}}}
\multiputlist(35,10)(20,0){$\scriptscriptstyle \alpha_1$,
$\scriptscriptstyle \alpha_2$,$\scriptscriptstyle \alpha_i$,
$\scriptscriptstyle \alpha_{n{-}1}$,$\scriptscriptstyle \alpha_n$}
\put(73,-15){$\scriptstyle \uparrow$}
\end{picture}
\begin{picture}(10,50)(-200,-30)
\put(0,0){$B_n$}
\multiput(35,0)(20,0){4}{\circle{8}}
\put(115.5,0){\circle*{8.5}}
\multiputlist(45,0)(20,0)
{{\line(1,0){12}},{\dashline{2}(-5,0)(5,0)},{\line(1,0){12}}}
\multiputlist(35,10)(20,0){$\scriptscriptstyle \alpha_1$,
$\scriptscriptstyle \alpha_2$,$\scriptscriptstyle \alpha_{n{-}2}$,
$\scriptscriptstyle \alpha_{n{-}1}$,$\scriptscriptstyle \alpha_n$}
\put(99,1.5){\line(1,0){13}}
\put(99,-1.5){\line(1,0){13}}
\put(32.5,-15){$\scriptstyle \uparrow$}
\end{picture}

\begin{picture}(10,50)(-20, -40)
\put(0,0){$C_n$}
\multiput(35,0)(20,0){4}{\circle*{8.5}}
\put(115.5,0){\circle{8}}
\multiputlist(45,0)(20,0)
{{\line(1,0){12}},{\dashline{2}(-5,0)(5,0)},{\line(1,0){12}}}
\multiputlist(35,10)(20,0){$\scriptscriptstyle \alpha_1$,
$\scriptscriptstyle \alpha_2$,$\scriptscriptstyle \alpha_{n{-}2}$,
$\scriptscriptstyle \alpha_{n{-}1}$,$\scriptscriptstyle \alpha_n$}
\put(99,1.5){\line(1,0){12.5}}
\put(98.5,-1.5){\line(1,0){13}}
\put(114,-15){$\scriptstyle \uparrow$}
\end{picture}
\begin{picture}(10,50)(-200,-40)
\put(0,0){$D_n$}
\multiput(35,0)(20,0){4}{\circle{8}}
\put(115.5,15){\circle{8.5}}
\put(115.5,-15){\circle{8.5}}
\multiputlist(45,0)(20,0)
{{\line(1,0){12}},{\dashline{2}(-5,0)(5,0)},{\line(1,0){12}}}
\multiputlist(35,10)(20,0){$\scriptscriptstyle \alpha_1$,
$\scriptscriptstyle \alpha_2$,$\scriptscriptstyle \alpha_{n{-}3}$,
$\scriptscriptstyle \alpha_{n{-}2}$,$\scriptscriptstyle $}
\put(107,25){$\scriptscriptstyle \alpha_{n{-}1}$}
\put(110,-25){$\scriptscriptstyle \alpha_{n}$}
\put(99,1.5){\line(4,3){13}}
\put(98.5,-1.5){\line(4,-3){13.5}}
\put(33.5,-15){$\scriptstyle \uparrow$}
\put(125,13){$\scriptstyle \leftarrow$}
\put(125,-17){$\scriptstyle \leftarrow$}
\end{picture}

\begin{picture}(10,40)(-20,-30)
\put(0,0){$E_6$}
\multiput(35,0)(20,0){5}{\circle{8}}
\multiputlist(45,0)(20,0)
{{\line(1,0){12}},{\line(1,0){12}},{\line(1,0){12}},{\line(1,0){12}}}
\multiputlist(35,10)(20,0){$\scriptscriptstyle \alpha_1$,
$\scriptscriptstyle \alpha_2$,$\hspace{11pt}\scriptscriptstyle \alpha_3 $,
$\scriptscriptstyle \alpha_5$,$\scriptscriptstyle \alpha_6$}
\put(70.5,30){$\scriptscriptstyle \alpha_4$}
\put(74.9,20){\circle{8}}
\put(74.5,4.1){\line(0,1){12.2}}
\put(33,-15){$\scriptstyle \uparrow$}
\put(113,-15){$\scriptstyle \uparrow$}
\end{picture}
\begin{picture}(10,40)(-200,-30)
\put(0,0){$E_7$}
\multiput(35,0)(20,0){6}{\circle{8}}
\multiputlist(45,0)(20,0)
{{\line(1,0){12}},{\line(1,0){12}},{\line(1,0){12}},{\line(1,0){12}},{\line(1,0){12}}}
\multiputlist(35,10)(20,0){$\scriptscriptstyle \alpha_1$,
$\scriptscriptstyle \alpha_2$,$\scriptscriptstyle \alpha_3 $,
$\hspace{8pt}\scriptscriptstyle \alpha_4$,$\scriptscriptstyle \alpha_6$,$\scriptscriptstyle \alpha_7$}
\put(70.5,30){$\hspace{18pt}\scriptscriptstyle \alpha_5$}
\put(94.2,20){\circle{8}}
\put(94.4,4.8){\line(0,1){12.2}}
\put(33,-15){$\scriptstyle \uparrow$}
\end{picture}
\caption {Dynkin diagrams with the non compact roots marked.}
\label{dynkin}
\end{figure}
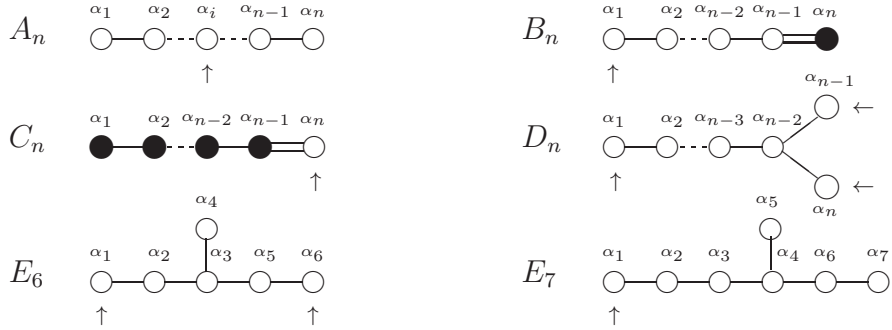

We mark in Figure \ref{dynkin} the non compact roots $\phi$ in the Dynkin diagram of the complexified Lie algebras $\g$.

\section{Nijenhuis tensor on compact hermitian symmetric spaces}
We collect here the basic facts that we need about Poisson-Nijenhuis geometry (see for instance \cite{Bonechi} for more details) and we introduce the examples we are going to study. For background in Poisson-Lie groups and Poisson homogeneous spaces see for instance \cite{Lu1}. Let $M$ be a smooth manifold, $\omega$ be a symplectic form and $\pi$ be a Poisson tensor on $M$ that are compatible, {\it i.e.}
\begin{equation}\label{PN_compatibility}
[\omega^{-1},\pi]=0\;.
\end{equation}
As a consequence the $(1,1)$-tensor $N=\pi\circ\omega$ has vanishing torsion
\begin{equation}\label{Nijenhuis_torsion}
 T(N)(v_1,v_2) = N([Nv_1,v_2]+[v_1,Nv_2]-N[v_1,v_2])-[Nv_1,Nv_2]=0\;.
\end{equation}

We remark that the PN geometry defined by (\ref{PN_compatibility}) has been recently
generalized in terms of compatibility between the Nijenhuis tensor and Dirac structures, together with the lift to Lie groupoids, in \cite{BDN}. 

Let $\Omega(M)$ denote the graded ring of differentiable forms and let $\iota_N:\Omega(M)\rightarrow\Omega(M)$ the degree zero derivation defined as $\iota_N(f)=0$ and $\iota_N(\nu)=N^t\nu$ for all $f\in C^\infty(M)$ and $\nu\in\Omega^1(M)$. As a consequence of the vanishing of the Nijenhuis torsion
\begin{equation}
\label{Nijenhuis_complex}
d_N =[d,\iota_N]=d\iota_N-\iota_N d
\end{equation}
squares to zero.  Remark that $d_Nf=Ndf$ for $f\in C^\infty(M)$. We refer to the complex $\Omega_N=(\Omega(M),d_N)$ as the {\it Nijenhuis complex}. We have clearly that $[d,d_N]=0$; Hamiltonian forms defined in (\ref{hamiltonian_forms}) are then those one forms that are closed with respect to both $d$ and $d_N$.

The symplectic form $\omega:TM\rightarrow T^*M$ defines an isomorphism of complexes between $(\Omega(M),d_N)$ and the complex $(\Gamma(\Lambda TM), d_\pi)$ computing the Lichnerowicz Poisson cohomology of $\pi$.

Let $K$ a compact simple Lie group; we denote with $\k$  its Lie algebra and with $\t$ a choice of Cartan subalgebra. Let $\phi$ be a non compact root, $\k_\phi\subset\k$ the corresponding Lie subalgebra and $K_\phi\subset K$ the subgroup. Let $\rho_\phi$ be the generator of the one dimensional center $\z(k_\phi)$ normalized by $\phi(\rho_\phi)=i$. The homogeneous space $K/K_\phi$ can be realized as the adjoint orbit $M_\phi$ of $\rho_\phi$, endowing it with the KKS symplectic form $\omega$. The moment map of the $K$-action is denoted 
$$\mu:M_\phi\rightarrow \k^*\equiv\k,$$
where the identification is done thanks to the Killing form. 

The standard Poisson-Lie structure $\pi_K$ on $K$ is defined
for each $k\in K$ as
$$\pi_K(k) = l_k(r)-r_k(r)$$
where $r$ is the classical $r$-matrix
$$
r= \frac{i}{2} \sum_{\alpha\in\Delta^+} e_\alpha\wedge e_{-\alpha}\;.
$$
The subgroup $K_\phi$ is a Poisson subgroup and induces a Poisson structure
on $M_\phi$ that we denote with $\pi$ and we call the {\it Bruhat-Poisson structure}. The homogeneous $K$-action on $M_\phi$ is a Poisson action, {\it i.e.} for each $X\in\k$ we have that $$[X^\sharp,\pi ]=\delta(X)^\sharp\;,$$
where $\delta:\k\rightarrow \Lambda^2\k$ is the cocycle encoding the dual Lie algebra on $\k^*$ (see Thm.2.6 in \cite{LuWe}).

The compatibility between the KKS symplectic from and the Bruhat-Poisson structure $\pi$ has been proven in \cite{KRR}.

We skip all details concerning the Bruhat-Poisson structure, that can be found for instance in \cite{LuWe}. In this paper we will only need the following formula that has been proven in Theorem 6.1 of \cite{BQT}
\begin{equation}
 \label{the_mother_of_all_formulas}
 d_N\mu = d\mu -[J(d\mu),\mu]\;.
\end{equation}

The following Lemma is a consequence of the fact that the $K$-action on $M_\phi$ is a Poisson action. Let $\k_1\subset\k$ any Lie subalgebra and let $(\Omega^{\k_1},d)$ be the complex of $\k_1$ invariant forms.

\begin{lemma}\label{invariant_forms}
$(\Omega^{\k_1},d_N)$ is a subcomplex of the Nijenhuis complex.
\end{lemma}
{\it Proof}. Since the $\k$-action (and then $\k_1$) is hamiltonian, the cochain map $\omega:(\Omega,d_N)\rightarrow (\Gamma(\Lambda TM),d_\pi)$ exchanges the $\k$-actions via Lie derivatives on forms and multivector fields. So it is enough to prove the statement on the LP complex. Indeed, for each $X\in\k_1$ and $A\in\Gamma(\Lambda TM)$ such that $L_{X^\sharp}(A)=0$ we have
$$
L_{X^\sharp}(d_\pi(A)) = [X^\sharp,[\pi,A]] = [[X^\sharp,\pi],A]=[\delta(X)^\sharp,A]=0\,.
$$
\qed

\section{Nijenhuis operator and invariant polynomials}
We saw in Lemma \ref{invariant_forms} that invariant forms with respect to $\k$ (and to any subalgebra) form a subcomplex of the Nijenhuis complex. Here we ask the same question about the subcomplex generated by invariant polynomials.

We limit ourselves to a class of subalgebras that is compatible with the complex structure $J$ defined in (\ref{complex_structure}). Let $\k_1\subset \k$ be a Lie subalgebra and $K_1\subset K$ be the subgroup integrating it.
It is clear that $\k_1^\perp$ is invariant under the adjoint action of $\k_1$. Let us ask now that both $\k_1$ and $\k_1^\perp$ are $J$-invariant and that
for each $X\in\k_1$ we have that
\begin{equation}
 \label{hypothesys_subalgebra}
 [\ad_X, J|_{k_1^\perp}]=0  \;.
\end{equation}
We denote with $\pr_{\k_1}$ and $\pr_{\k_1^\perp}$ the orthogonal projections to $\k_1$ and $\k_1^\perp$ respectively.

\begin{lemma}\label{J_derivation}
If $\k_1$ satisfies the conditions (\ref{hypothesys_subalgebra}) then, for each $\xi,\eta\in\k_1^\perp$ we have that
$$
\pr_{\k_1}([J\xi,\eta]+[\xi,J\eta])=0
$$
\end{lemma}
{\it Proof}. For each $X\in\k_1$ we have that
$$
( X, [J\xi,\eta]+[\xi,J\eta]) = ( [X,J\xi],\eta) +( [X,\xi],J\eta) =  ( J[X,\xi],\eta ) +( [X,\xi],J\eta)=0
$$
where in the second equality we used (\ref{hypothesys_subalgebra}) and in the last one the antisymmetry of $J$ with respect to the Killing form. \qed

\smallskip
We denote with $\mu_{\k_1}=M_\phi\rightarrow\k_1$ the moment map and let $\mu=\mu_{\k_1}+\mu_{\k_1^\perp}$ the decomposition of $\k$-moment map given by $\k=\k_1\dot + \k_1^\perp$. It is clear that
for each $X\in\k_1$
\begin{equation}
 \label{invariant_operator1}
 X^\sharp(\mu_{\k_1^\perp}) = [ X,\mu_{\k_1^\perp}]\,,\,\,\, X^\sharp(J\mu_{\k_1^\perp}) = [ X,J\mu_{\k_1^\perp}]\,,
\end{equation}
where the second equality is a consequence of (\ref{hypothesys_subalgebra}).

Let us consider the space $\Symm^r(\k_1)^{\k_1}$ of invariant polynomials that are homogeneous of degree $r$ and let $\Symm(\k_1)^{\k_1}=\oplus_{r\geq 0}\Symm^r(\k_1)^{\k_1}\subset C^\infty(\k^*)$.
If $p\in\Symm^r(\k_1)^{\k_1}$, we denote with $\tilde{p}$ the corresponding $r$-linear application $\k_1\otimes\ldots\otimes\k_1\rightarrow\C$.
In order to avoid cumbersome notations, we use the same symbol $p$ to denote $\mu^*_{\k_1}(p)\in C^\infty(M_\phi)$. We denote with $(\Omega^{\k_1}_{b},d)\subset (\Omega^{\k_1},d)$ the subcomplex of $\k_1$-basic forms and with $(\Omega^{\k_1}_{pol},d)\subset (\Omega^{\k_1}_b,d)$ the subcomplex generated by $S(\k_1)^{\k_1}$.
We would like to study the conditions under which $\Omega^{\k_1}_{pol}$ is a subcomplex of the Nijenhuis complex. We prove the following preliminary result.

\begin{prop}\label{thm_basic_forms}
\begin{itemize}
\item[$i$)] If $\k_1$ satisfies conditions (\ref{hypothesys_subalgebra}) then
for each $p\in S^r(\k_1)^{\k_1}$ we compute
\begin{equation}\label{action_Nijenhuis}
d_Np= dp - r \tilde{p}(\mu_{\k_1},\dots,\mu_{\k_1},dA_{\k_1})\,,
\end{equation}
where we defined
\begin{equation}\label{fundamental_object}
A_{\k_1} = \frac{1}{2}\pr_{\k_1}[J\mu_{\k_1^\perp},\mu_{\k_1^\perp}]\,.
\end{equation}
As a consequence $d_Np$ is a $\k_1$-invariant
form.
\item[$ii$)] If
\begin{equation}\label{sufficient_condition}
 [A_{\k_1},\mu_{\k_1}]=0
\end{equation}
then $d_Np$ is a $\k_1$-basic form.
\end{itemize}
\end{prop}
{\it Proof.} We compute from (\ref{the_mother_of_all_formulas})
$$
d_N\mu_{\k_1} = d\mu_{\k_1} -\pr_{\k_1}[Jd\mu,\mu] = d\mu_{\k_1} -[Jd\mu_{\k_1},\mu_{\k_1}] - \pr_{\k_1}[Jd\mu_{\k_1^\perp},\mu_{\k_1^\perp}]\;.
$$
From Lemma \ref{J_derivation} we see that
$$
\pr_{\k_1}[Jd\mu_{\k_1^\perp},\mu_{\k_1^\perp}]= d A_{\k_1}\;.
$$

As a consequence of (\ref{invariant_operator1}) we have that for each $X\in\k_1$
\begin{equation}
 \label{invariant_operator2}
 X^\sharp(A_{\k_1}) = [ X,A_{\k_1}]\,.
\end{equation}
Let $p\in\Symm(\k_1)^{\k_1}$. We compute
\begin{eqnarray*}
d_Np &=&\sum_{i=1}^r \tilde{p}(\mu_{\k_1},\dots, d_N\mu_{\k_1},\ldots,\mu_{\k_1})=
r \tilde{p}(\mu_{\k_1},\dots,\mu_{\k_1}, d_N\mu_{\k_1})\cr
&=&
dp - r \tilde{p}(\mu_{\k_1},\dots,\mu_{\k_1},[Jd\mu_{\k_1},\mu_{\k_1}] + dA_{\k_1}) \;.
\end{eqnarray*}
By using the invariance of $\tilde{p}$ we get (\ref{action_Nijenhuis}).

Since $A_{\k_1}$ satisfies (\ref{invariant_operator2}) and $p$ is an invariant polynomial, it follows that for each $X\in\k_1$
$$
L_{X^\sharp}(d_Np) = 0\:,
$$
so that $d_Np$ is invariant.
Let us assume that $[\mu_{\k_1},A_{\k_1}]=0$; from (\ref{invariant_operator2}) we now compute
$$\iota_{X^\sharp}(d_Np) = - r \tilde{p}(\mu_{\k_1},\dots,\mu_{\k_1},[X,A_{\k_1}])
= - r\sum_{i=1}^{r-1} \tilde{p}(\mu_{\k_1},\ldots,[A_{\k_1},\mu_{\k_1}],\ldots, X)=0 .
$$\qed

\begin{cor}
 Let $q,p\in S(\k_1)^{\k_1}$ then
 \begin{equation}\label{commutation}
 \{q,p\}_{\pi} =0\;.
 \end{equation}
Let us suppose that (\ref{sufficient_condition}) holds and let
$\k_2\subset\k_1$. Then (\ref{commutation}) holds for each $q\in S(\k_2)^{\k_2}$ and $p\in S(\k_1)^{\k_1}$.

\end{cor}
{\it Proof}. Let $v_q=\omega^{-1}(dq)$ be the hamiltonian vector field of $q$. We compute from (\ref{action_Nijenhuis})
\begin{eqnarray}
\{q,p\}_{\pi}&=&\iota_{v_q}(d_Np) =- r \tilde{p}(\mu_{\k_1},\ldots,\mu_{\k_1},v_{q}(A_{\k_1})) =-r \tilde{p}(\mu_{\k_1},\ldots,\mu_{\k_1},[X_q,A_{\k_1}]) \cr
&=& r \sum_{i=1}^{r-1} \tilde{p}(\mu_{\k_1},\ldots [X_q,\mu_{\k_1}],\ldots A_{\k_1}) \label{commutation_first_step}
\end{eqnarray}
where $X_q\in\k_1\subset\k$ is defined as
$$
(Y, X_q ) = \tilde{q}(\mu_{\k_1},\ldots,\mu_{\k_1},Y)
$$
for each $Y\in\k_1$ so that $X^\sharp_q=v_q$ and the second step follows from (\ref{invariant_operator2}). It is now clear that
$[X_q,\mu_{\k_1}]=0$, since for each $Y\in\k_1$ we have that
$$ ( Y,[\mu_{\k_1},X_q])= ( [Y,\mu_{\k_1}],X_q) =
\tilde{q}(\mu_{\k_1},\ldots,\mu_{\k_1},[Y,\mu_{\k_1}])=0$$
since $q$ is invariant.

If (\ref{sufficient_condition}) holds, then by using the invariance of $\tilde{p}$ in the last equality of the first line of (\ref{commutation_first_step}), we get that (\ref{commutation}) holds also when $q\in S(\k_2)^{\k_2}$ and $p\in S(\k_1)^{\k_1}$. \qed

\bigskip
\bigskip

\section{The Thimm chain of subalgebras}\label{thimm_section}

We recall here the Thimm method for constructing integrable models (see \cite{GS1} for details). Let us consider $\k\supset\k_1\supset\k_2$ nested Lie subalgebras and let $q\in S(\k_1)^{\k_1}$ and $p\in S(\k_2)^{\k_2}$. Since $\mu_{\k_2}=p_{21}\circ\mu_{\k_1}$, where $p_{21}:\k_1^*\rightarrow\k_2^*$ we have that
$$
\{\mu^*_{\k_2}p,\mu^*_{\k_1}q\} = \{\mu^*_{\k_1}p_{21}^*p,\mu^*_{\k_1}q\} =
\mu_{\k_1}^*\{p_{21}^*p,q\}_{\k_1^*}=0
$$
since $S^{\k_1}(\k_1)\subset C^\infty(\k_1^*)$ is the centre of the Poisson algebra. Let us consider now a chain of Lie subalgebras
$$
\k\equiv \k_0\supset \k_1\supset \ldots\k_n\;
$$
and let $K_i\subset K$ the corresponding subgroups. Let us denote with
\begin{equation}\label{thimm_hamiltonians}S(\k\supset\k_1\ldots\supset \k_n)= \cup_i S(\k_i)^{\k_i}\end{equation}
the union of invariant polynomials of any subalgebra of the chain. By the above observation, all these polynomials are in involution.

Let us discuss if they are independent and so form an integrable model. The following two properties are equivalent: $i$) the ring of invariant functions $C^\infty(M)^{\k_1}$ is abelian with respect to the Poisson bracket; $ii$) for each coadjoint orbit $O\subset\k^*$ the symplectic reduction $\mu_\k^{-1}(O)/K$ is a point. If one of the above condition (and then both) is true then the $K$-action is said to be {\it multiplicity free} (\cite{GS}). If for each $K_i$- orbit $O\subset M$, the action of $K_{i-1}$ on $O$ is multiplicitely free then the hamiltonians in (\ref{thimm_hamiltonians}) define an integrable model.

Let us consider now the case $M=M_\phi$; we want to see the conditions under which the hamiltonians (\ref{thimm_hamiltonians}) commute also with respect to the Poisson structure $\pi$. Let us choose a non compact root $\phi_1$ of $\k$ and let $\rho_1$ the normalized generator of the one dimensional centre $\z_1$ of $\k_{\phi_1}$.
In general this root $\phi_1$ can be different from $\phi$ involved in the definition of the hermitian symmetric space.

Let $\k_1=\k_{\phi_1}= \k_{\phi_1}'\oplus \z_1$ with $\k_{\phi_1}'$ being simple. Let $\phi_2$ be a non compact root of $\k_{\phi_1}'$ with $\rho_2$ the normalized generator of the one dimensional center $\z_2$ of $\k_{\phi_2}\subset\k_{\phi_1}'$. Let $\k_2 = \k_{\phi_2}\oplus\z_1=\k_{\phi_2}'\oplus\z_2\oplus\z_1$. We get the decomposition
$$
\k = \k_2\dot{+}\k_2^\perp=\k_2 \dot{+} (\k_2^\perp\cap\k_1) \dot{+} \k_1^\perp=\k_2\dot{+}(\k_2^\perp\cap \k_{\phi_1}')\dot{+}\k_1^\perp\;,
$$
where the third equality follows because we included $\z_1$ in the definition of $\k_2$. We remark that the perpendicular is always taken in $\k$.
We can iterate the procedure and choose $\phi_i$ non compact root of $\k_{\phi_{i-1}}'$ with $\rho_i$ the normalized generator of the centre $\z_i$ of
$\k_{\phi_i}\subset \k_{\phi_{i-1}'}$ and define
\begin{equation}\label{generic_subalgebra_chain}
\k_i = \k_{\phi_i}\oplus\z_{i-1}\oplus\ldots \oplus\z_1\;.
\end{equation}
We get the $i-th$ decomposition of $\k$
\begin{equation}
 \label{ith_decomposition}
 \k = \k_i \dot{+} (\k_i^\perp\cap \k_{\phi_{i-1}}') \dot{+} \k_{i-1}^\perp = \k_i \dot{+} (\k_i^\perp\cap \k_{\phi_{i-1}}') \dot{+} (\k_{i-1}^\perp\cap \k_{\phi_{i-2}}') \ldots \dot{+}\k_1^\perp\;.
\end{equation}
We remark that $J$ acts on each addend $\k_i^\perp\cap \k_{\phi_{i-1}'}$ of (\ref{ith_decomposition}) as $\ad_{\rho_{i}}$.

\begin{lemma}
Every subalgebra $\k_i$ satisfies condition (\ref{hypothesys_subalgebra}).
\end{lemma}
{\it Proof}.  By construction the Cartan subalgebra $\t$ is included in $\k_i$ for each $i$. Moreover, given any root $\alpha$, either the root vector $e_\alpha\in\k_i$ or $e_\alpha\in\k_i^\perp$. This implies that $\k_i$ and $\k_i^\perp$ are $J$-invariant.  Equation (\ref{hypothesys_subalgebra}) is straightforward if $X\in\t\subset\k_i$. Let $X=e_\alpha\in\k_i$ and $e_\beta\in \k_i^\perp$ be two root vectors.
If $\alpha+\beta$ is not a root then both terms are zero; if $\alpha+\beta$ is a root then it is non compact and it will be positive if and only if $\beta$ is positive so that
$[X,J(e_{\beta})]=J[X, e_\beta]$. \qed

\smallskip
We will see in the examples that, for each simple compact Lie algebra $\k$, it is possible to find such a chain that ends with $\k_n=\t$ the Cartan subalgebra of $\k$. Such a choice will not be unique. We call such a chain
\begin{equation}\label{thimm_chain}
 \k\supset \k_1\supset \k_2\ldots\supset \k_n=\t
\end{equation}
a Thimm chain.

\section{The $\phi\,$-minimal representations}

Let $\phi$ be a non compact root of the simple Lie algebra $\k$ and let $\rho_\phi$ the normalized generator of the center $\z$ of $\k_\phi$. Let $R_\Lambda:\k\rightarrow {\rm End}(V(\Lambda))$ be an irreducible representation of highest weight $\Lambda\in\t^*$. It is easy to see that $V(\Lambda)$ decomposes in eigenspaces of $\rho_\phi$ as
\begin{equation}\label{general_decomposition_irrep}
V(\Lambda) = \oplus_{\ell\geq 0} V(\Lambda)_\ell
\end{equation}
where $V(\Lambda)_\ell$ are $\k_\phi$-invariant, $\rho_\phi|_{V(\Lambda)_\ell}= \Lambda_\phi-i\, \ell$ with $\Lambda_\phi\equiv \Lambda(\rho_\phi)$ and $e_\alpha: V(\Lambda)_\ell\rightarrow V(\Lambda)_{\ell- 1}$ for every non compact positive root $\alpha$.

\begin{defn}\label{phi_minimal_rep}
We say that $R_\Lambda$ is $\phi\,$-{\it minimal} if the sum in (\ref{general_decomposition_irrep}) runs over $\ell=0,1$ and $V(\Lambda)_0$ is irreducible. In this case we denote $V_+ \equiv V(\Lambda)_0$ and $V_-\equiv V(\Lambda)_1$.
\end{defn}

Let now $\phi$ be the noncompact root defining the Hermitian symmetric space $M_\phi$ and let us suppose that $V(\Lambda)$ is $\phi\,$-minimal. Since $M_\phi$ is the adjoint orbit of $\rho_\phi$, the moment map $\mu_\Lambda=R_\Lambda(\mu)$ in the representation 
$R_\Lambda$ satisfies then \begin{equation}\label{quadratic_relation_moment_map}
(\mu_\Lambda-\Lambda_\phi)(\mu_\Lambda-\Lambda_\phi+i)=\mu_\Lambda^2 -(2\Lambda_\phi-i)\mu_\Lambda +\Lambda_\phi (\Lambda_\phi-i)=0\;.
\end{equation}

Let us consider now the Thimm chain (\ref{thimm_chain}) and let us suppose that $V(\Lambda)$ is also $\phi_1$-minimal
as $V(\Lambda)=V_+^1\oplus V_-^1$. We can go on along the Thimm chain (\ref{thimm_chain}) and let us assume that $V_+^1$ is $\phi_2$-minimal as
$V_+^1= V_+^2\oplus V_-^2$ and so on. At level $i$, we assume that $V_+^{i-1}$ is $\phi_i$-minimal so that $V^{i-1}_+ = V^{i}_+\oplus V^{i}_-$. At each level $i$, we get the decomposition of $V$ into $\k_i$-representation as
\begin{equation}\label{V_decomposition}
V(\Lambda)= W^i_+\oplus W^i_-
\end{equation}
where $W^i_+= V^i_+$ and $W^i_-= V^i_-\oplus V^{i-1}_-\ldots \oplus V^{1}_-$. If these conditions are satisfied we say that the representation $R_\Lambda$ is minimal with respect to the chain (\ref{thimm_chain}).

\smallskip
\begin{lemma}\label{decomposable_representation_lemma}
If the representation $R_\Lambda$ is minimal with respect to the Thimm chain (\ref{thimm_chain}) then each $\xi\in\k_i^\perp$ decomposes with respect to (\ref{V_decomposition}) as
$$
R_\Lambda(\xi) = \left(\begin{array}{cc}
                             0 & \nu\cr
                             -\nu^\dagger & \lambda
                            \end{array}
\right)\,,\;\;\;
R_\Lambda(J(\xi)) = \left(\begin{array}{cc}
                             0 & i\nu\cr
                             i\nu^\dagger & J(\lambda)
                            \end{array}
\right)
$$
\end{lemma}
{\it Proof}. 
Let us consider the decomposition (\ref{V_decomposition}) of $V$ in $\k_i$ representations. It is clear that $\k_i^\perp:W_+^i\rightarrow W_-^i$. Indeed
if $\xi\in \k_i^\perp\cap\k_{\phi_{i-1}}'$ then $\xi: W_+^i=V_+^i\rightarrow V_-^i\subset W^i_-$; if $\xi \in\k_{i-1}^\perp$ then $\xi: W_+^i=V_+^i\subset V_+^{i-1}\rightarrow V_-^{i-1}\subset W^i_-$. Let $\alpha$ be a positive root such that $e_\alpha\in (\k_i^\perp\cap\k_{\phi_{i-1}}')\otimes\C$, then $e_\alpha|_{W_+^i} = e_\alpha|_{V_+^i}=0$; if $e_\alpha\in\k_{i-1}^\perp\otimes\C$ then $e_\alpha|_{W_+^i}=e_\alpha|_{V_+^i}= (e_\alpha|_{V^{i-1}_+})|_{V_+^i}=0$. \qed

\medskip

Let us consider now the the ring of $\k_i$-invariant polynomials generated by the traces in the $W^i_+$ representation: let us define
\begin{equation}\label{trace_polinomials}
I_r^{(i)}(X) = \frac{i^r}{r}\Tr_{W^i_+}(X^r)\;.
\end{equation}

\begin{thm}\label{explicit_formula}
We compute
\begin{equation}
\label{action_Nijenhuis_indecomposable}
d_NI_r^{(i)} = -2i\Lambda_\phi d I_r^{(i)} + 2dI_{r+1}^{(i)}\;.
\end{equation}
\end{thm}
{\it Proof}. We then have to compute (\ref{fundamental_object}) in the representation $W^i_+$. If we write
$$
R_\Lambda(\mu_{\k_i}) = \left(\begin{array}{cc}
                             R_{W^i_+}(\mu_{\k_i}) & \nu_i\cr
                             -\nu_i^\dagger & R_{W^i_-}(\mu_{\k_i}) +\lambda_i
                            \end{array}
\right)\,,\;\;\;
$$
then by using Lemma \ref{decomposable_representation_lemma} we can compute $A_{\k_i}$ defined in (\ref{fundamental_object}) in the representation $W_+^i$ as
$$
R_{W^i_+}(A_{\k_i})=-i\nu_i\nu_i^\dagger\;.
$$
From (\ref{quadratic_relation_moment_map}) we find that
\begin{equation}\label{fundamental_from_quadratic}
R_{W^i_+}(A_{\k_i}) = -i R_{W^i_+}(\mu_{\k_i})^2 +(2i\Lambda_\phi+1)R_{W^i_+}(\mu_{\k_i}) -\Lambda_\phi(1+i\Lambda_\phi)\;.
\end{equation}
Let us now insert it in (\ref{action_Nijenhuis}) with $p=I_r$ and get
\begin{eqnarray*}
d_NI_r^{(i)} &=& dI_r^{(i)} + r\tilde{I_r}^{(i)}\left(\mu_{\k_i},\ldots, id(\mu_{\k_i}^2)-(2i\Lambda_\phi+1)d\mu_{\k_i}\right)\cr
&=& dI_r^{(i)} + i^{r+1}\Tr_{W^i_+}(\mu_{\k_i}^{r-1}d(\mu_{\k_i}^2))-(2i\Lambda_\phi+1) dI_r^{(i)}
\;,
\end{eqnarray*}
so that we finally get (\ref{action_Nijenhuis_indecomposable}). \qed

\medskip
From (\ref{action_Nijenhuis_indecomposable}) it is now easy to prove the following corollary. 

\begin{cor}\label{cor_indecomposable}
\begin{itemize}\item[$i$)] The subcomplex $(\Omega_{\k_i}^{pol},d)$ of the de Rham complex generated by $S(\k_i)^{\k_i}$, defines also a
subcomplex $(\Omega_{\k_i},d_N)$ of the Nijenhuis complex (\ref{Nijenhuis_complex});
\item[$ii$)] if $\lambda$ is an eigenvalue of the moment map $\mu_{k_i}$ then
\begin{equation}\label{Nijenhuis_eigenvalue}
\tilde{\lambda} = 2i(\lambda - \Lambda_\phi)
\end{equation}
is a Nijenhuis eigenvalue.
\end{itemize}
\end{cor}
In the next section we will prove that such minimal representations exist in the classical symmetric spaces so that the above Corollary reproduces the diagonalization proved case by case in \cite{BQT}.

\section{Existence of $\phi\,$-minimal representations}
Here we discuss the existence of $\phi\,$-minimal representations. We give the explicit examples for the classical groups and prove that they do not exist in the exceptional cases.

\subsection{AIII}
Let $M_\phi=SU(n)/S(U(k)\times U(n-k))$ and let us consider the root ordering as in $A_{n-1}$ Dynkin diagram in Figure \ref{dynkin}; the non compact root is $\phi=\alpha_k$. In the fundamental representation $f_{\mathfrak{su}(n)}$ of $\mathfrak{su}(n)$,
\begin{equation}\label{rho_AIII}
    \rho_{\alpha_k} = \frac{i}{n}\begin{pmatrix}
    (n-k)1_k &0\\
    0 & -k\,1 _{n-k}
    \end{pmatrix},
\end{equation}
so that $f_{\su(n)}$ is $\alpha_k$-minimal and  decomposes in representations of $\s(\u(k)\oplus\u(n-k))$ as
$f_{\mathfrak{su}(n)} = V_+ \oplus  V_-$, where
\begin{equation*}
    \begin{aligned}
    &  V_+ = (f_{\mathfrak{u}(k)},0_{\mathfrak{u}(n-k)}),\\
    & V_- =  (0_{\mathfrak{u}(k)}, f_{\mathfrak{u}(n-k)}),
    \end{aligned}
\end{equation*}
with $f_{\mathfrak{u}(k)}$ and $0_{\mathfrak{u}(k)}$ respectively the fundamental and the trivial representation. The eigenvalue of $\rho_{\alpha_k}$ on $V_+$ is then
\begin{equation}\label{lambda_AIII}
 \Lambda_{\alpha_k}=i\frac{n-k}{n}
\end{equation}

Now, let us consider the Thimm chain of subalgebras
\begin{equation}\label{thimm_chain_a}
    \mathfrak{s}\u(n)\supset \mathfrak{s}(\u(n-1)\oplus \u(1))\supset ...\supset \mathfrak{s}(\u (n-i)\oplus \u (1)^i) \supset... \supset \u(1)^{n-1}.
\end{equation}
Let us choose $\phi_1=\alpha_{n-1}$ so that $f_{\su(n)}$ decomposes with respect to $\s(\u(n-1)\oplus\u(1))$ as $f_{\mathfrak{su}(n)} = V^1_+ \oplus  V^1_-$, with
\begin{equation*}
 V^1_+ = \left( f_{\mathfrak{u}(n-1)},0\right).
\end{equation*}
We obtain for the $i$-th step of the chain (\ref{thimm_chain_a}) that $f_{\su(n)}$ decomposes with respect to
$\s(\u(n-i)\oplus\u(1)^i)$ as
\begin{equation*}
    f_{\su(n)} = W^{i}_+ \oplus W^{i}_-,
\end{equation*}
with
\begin{equation*}
    W^{i}_+ \equiv \Big( f_{\mathfrak{u}(n-i)},\underbrace{0,...,0}_{i-\text{times}}\Big),
\end{equation*}
that is $\phi_{i+1}=\alpha_{n-1-i}\,$-minimal, where $\alpha_{n-1-i}$ is a root of $\su(n-i)$. Then, the fundamental representation $f_{\mathfrak{su}(n)}$ of $\mathfrak{su}(n)$ is minimal with respect to the chain of subalgebras (\ref{thimm_chain_a}). Formula (\ref{action_Nijenhuis_indecomposable}) then gives
\begin{equation}
    d_NI_r^{(i)} = 2\frac{n-k}{n} dI_r^{(i)} +2 dI_{r+1}^{(i)}.
\end{equation}

\subsection{BDI}
Let us consider $M_\phi=SO(n+2)/SO(n)\times SO(2)$ and the root ordering of Dynkin diagrams $D_m$ and $B_m$ in Figure \ref{dynkin}, for $n+2=2m$ and $n+2=2m+1$ cases respectively. In these Dynkin diagrams, the non compact root defining the compact hermitian symmetric space is $\phi \equiv \alpha_1$. The fundamental representation is realized by anti-symmetric real matrices acting on $\mathbb{R}^{n+2}$. We can select a Cartan subalgebra as spanned by
\begin{equation*}
\begin{aligned}
&n+2=2m+1: \quad
    &&\mathfrak{t}= \left\lbrace \begin{pmatrix} 0 & 0 \\ 0 & a\otimes \sigma
    \end{pmatrix}, a=\mathrm{diag}(a_1,...,a_m), a_i\in \mathbb{R} \right\rbrace  \\
 &n+2=2m: \quad
    &&\mathfrak{t}= \left\lbrace a\otimes \sigma, \;a=\mathrm{diag}(a_1,...,a_m), a_i\in \mathbb{R}\right\rbrace,
\end{aligned}
\end{equation*}
with $\sigma = \begin{pmatrix} 0 & 1 \\ -1 & 0
    \end{pmatrix}$. Then, in both cases,
\begin{equation}
    \rho_\phi = \begin{pmatrix}
    0_{n} & 0 \\ 0 & \sigma
    \end{pmatrix},
\end{equation}
which clearly has 3 different eingenvalues. So, we can conclude that the fundamental representation is not $\phi\,$-minimal.

We are going to show that the spin representation is $\phi\,$-minimal. Let us consider the Clifford algebra $Cl(n+2,\mathbb{R})$, let $z_i$ coordinates on $\C^m$. The gamma matrices act on  $S^{(n+2)}=\wedge \C^m$  as
\begin{equation*}
    \Gamma_i = d\bar{z}_i\wedge, \qquad \Gamma_{\bar{i}}=\imath_{\partial \bar{z}_i}
\end{equation*}
together with $\Gamma_0=(-1)^{\rm deg}$ if $n$ is odd.
Then, the $\mathfrak{so}(n+2)$-spin representation is realized as
\begin{equation*}
    S(X) = \frac{1}{8}X_{ij}[\Gamma_i,\Gamma_j], \qquad X\in \mathfrak{so}(n+2).
\end{equation*}
In particular,
\begin{equation}
    S(\rho_\phi) = i\left( \Gamma_{\bar{m}}\Gamma_m -\frac{1}{2}\right).
\end{equation}
It has two different eingenvalues $\pm i/2$ so that  the spin representation is $\phi\,$-minimal, i.e. $S^{(n+2)}=V_+\oplus V_-$, where
\begin{equation*}
    V_+=\wedge \langle d\bar{z}_i, i=1,...,m-1\rangle, \qquad V_- = V_+(\Lambda)\otimes d\bar{z}_m,
\end{equation*}
which are $(S^{(n)},\pm i/2)$ representation of $\mathfrak{so}(n)\oplus \mathfrak{so}(2)$ respectively. Easily, one obtains that
\begin{equation}
    \Lambda_\phi = \frac{i}{2}.
\end{equation}

Now, let us consider the Thimm chains
\begin{equation}\label{thimm_chain_b}
    \begin{aligned}
    & \mathfrak{so}(2m) \supset...\supset \so(2m-2j)\oplus \so(2)^j \supset ...\supset \mathfrak{so}(2)^m,\\
    & \mathfrak{so}(2m+1) \supset ...\supset \so(2m+1-2j)\oplus \so(2)^j \supset...\supset\mathfrak{so}(3)\oplus \mathfrak{so}(2)^{m-1}.
    \end{aligned}
\end{equation}
This calculation can be iterated step-by-step along these chains and we get
\begin{equation*}
    W_+^j=\Big( S^{(n+2-2j)},\underbrace{\frac{i}{2},...,\frac{i}{2}}_{j-\text{times}}\Big),
\end{equation*}
with $S^{(n+2-2j)}$ the spin representation of $\mathfrak{so}(n+2-2j)$. Moreover, $ W_+^j(\Lambda)$ is $\phi_{j+1}=\alpha_{1}$-minimal, where $\alpha_1$ is the non compact root of $\so(n+2-2j)$. The spin representation is then minimal with respect to the chain in (\ref{thimm_chain_b}). Then, the Nijenhuis tensor acts on $\so(2n+2-2j)$-invariant polynomials $ I_r^{(j)}$, as defined in (\ref{trace_polinomials}), as
\begin{equation}
    d_NI_r^{(j)}= dI_r^{(j)}+2dI_{r+1}^{(j)}.
\end{equation}

\subsection{DIII}
Let us consider $M_\phi=SO(2n)/U(n)$ and the same root ordering of the $D_n$ case in Figure \ref{dynkin}. The non compact simple root defining this compact hermitian symmetric space is $\phi \equiv \alpha_n$. In the fundamental representation, the Lie subalgebra $\mathfrak{u}(n)$ is embedded in $\mathfrak{so}(2n)$ as
\begin{equation*}
    \mathfrak{u}(n) \ni A+iB \; \longrightarrow \; \begin{pmatrix} A & B \\ -B & A\end{pmatrix},
\end{equation*}
with $A,B  \in M_n (\mathbb{R})$ such those $A=-A^t,\, B=B^t$. The Cartan subalgebra can be chosen as
\begin{equation*}
    \mathfrak{t}= \left\lbrace \begin{pmatrix} 0_n & a \\ -a & 0_n    \end{pmatrix}, a = \mathrm{diag}(a_1,...,a_n), a_i\in \mathbb{R}\right\rbrace.
\end{equation*}
So,
\begin{equation}
    \rho_\phi= \begin{pmatrix} 0 & \frac{1}{2} 1_n \\ -\frac{1}{2} 1_n &0 \end{pmatrix}.
\end{equation}
Then, the fundamental representation $V=\C^{2n}$ is $\phi\,$-minimal, i.e. $V=V_+\oplus V_-$, where
\begin{equation*}
    V_+=\left\langle (a, +ia), a \in \mathbb{C}^n\right\rangle=f_{\mathfrak{u}(n)},
\end{equation*}
and, clearly,
\begin{equation}
    \Lambda_\phi=\frac{i}{2}.
\end{equation}

Moreover, if one considers the Thimm chain
\begin{equation}\label{thimm_chain_d}
    \mathfrak{so}(2n)\supset \mathfrak{u}(n)\supset...\supset \mathfrak{u}(n-i)\oplus \mathfrak{u}(1)^i\supset...\supset\mathfrak{u}(1)^{n},
\end{equation}
then, repeating the analysis for the case AIII one obtains that $V$ decomposes in representations of $\u(n-i)\oplus\u(1)^i$ as $V= W_+^i\oplus W_-^i$, where
\begin{equation*}
W^{i}_+ \equiv \Big( f_{\mathfrak{u}(n-i)},\underbrace{0,...,0}_{(i-1)-\text{times}}\Big),
\end{equation*}
and $W_+^i$ is $\phi_i=\alpha_{n-i-1}$-minimal, where $\alpha_{n-i-1}$ is the root of $\su(n-i)$. Then the fundamental representation is minimal with respect to the whole Thimm chain (\ref{thimm_chain_d}) and one obtains that
\begin{equation}
    d_NI_r^{(j)}= dI_r^{(j)}+2dI_{r+1}^{(j)}.
\end{equation}

\subsection{CI}
Let us consider $M_\phi=Sp(2n)/U(n)$ with the root ordering as in the $C_n$ case in Figure \ref{dynkin}. The non compact simple root defining the compact hermitian symmetric space CI is $\phi \equiv \alpha_n$. In the fundamental representation,
\begin{equation*}
  \mathfrak{sp}(2n) = \left\lbrace X(A,B)=\begin{pmatrix} A &B \\ -B^t & -A^t \end{pmatrix}, A,B\in M_n(\mathbb{C}), A=-A^t, B=B^t \right\rbrace,
\end{equation*}
while $\mathfrak{u}(n) = \left\lbrace X(A,0)\in \mathfrak{sp}(2n) \right\rbrace$. The Cartan subalgebra is spanned by matrices in $\mathfrak{u}(n)$ with $A$ diagonal. Then,
\begin{equation}
    \rho_\phi= \begin{pmatrix}  \frac{1}{2} 1_n & 0  \\ 0 & -\frac{1}{2} 1_n \end{pmatrix}.
\end{equation}
Then, the fundamental representation $f_{\sp(2n)}$ is $\phi\,$-minimal, i.e. $f_{\sp(2n)}=V_+\oplus V_-$, where
\begin{equation*}
    V_+=f_{\mathfrak{u}(n)},
\end{equation*}
and the eigenvalue of $\rho_\phi$ is
\begin{equation}
    \Lambda_\phi=\frac{i}{2}.
\end{equation}

Now, let us consider the Thimm chain
\begin{equation}\label{thimm_chain_c}
    \mathfrak{sp}(2n)\supset \mathfrak{u}(n)\supset\mathfrak{u}(n-1)\oplus\mathfrak{u}(1)\supset...\supset\mathfrak{u}(1)^n.
\end{equation}
At step $i+1$, $f_{\sp(2n)}$ decomposes in $\k_{i+1}=\u(n-i)\oplus\u(1)^i$-representations as $f_{\sp(2n)} = W_+^i\oplus W_-^i$,  where
\begin{equation*}
   W^{i}_+ \equiv \Big( f_{\mathfrak{u}(n-i)},\underbrace{0,...,0}_{i-\text{times}}\Big),
\end{equation*}
so that $W_+^i$ is $\alpha_{n-i-1}$-minimal, where $\alpha_{n-i-1}$ is a root of $\su(n-i)$. Then, the fundamental representation $f_{\sp(2n)}$  is minimal with respect to the Thimm chain in (\ref{thimm_chain_c}). Finally, one obtains that
\begin{equation}
    d_NI_r^{(j)}= dI_r^{(j)}+2dI_{r+1}^{(j)}.
\end{equation}

\subsection{EIII and EVII}
\begin{prop}\label{nogo_min_rep_exceptional}
There are no $\phi\,$-minimal representations of $e_6$ and $e_7$.
\end{prop}
{\it Proof}. Let $\Lambda$ be a dominant weight and let $V(\Lambda)$ be the finite dimensional irreducible representation of $\g$ with highest weight $\Lambda$ and highest weight vector $v_\Lambda$. If $V(\Lambda)$ is $\phi\,$-minimal then
$e_{-\beta}e_{-\alpha} v_\Lambda=0$ for each $\alpha,\beta\in\Delta_{nc}^+$ (otherwise the spectrum of $\rho_\phi$ contains the three distinct eigenvalues $\Lambda_\phi=\Lambda(\rho_\phi), \Lambda_\phi-i,\Lambda_\phi-2i$). Taking $\alpha=\beta$ this means that the string $\{\Lambda+j\alpha\}$ of weights extends at most to $-1\leq j\leq 0$, {\it i.e.}
\begin{equation}\label{first_condition}
(\Lambda,\alpha) \in\{0,1\}\;\;\;\forall\ \alpha\in\Delta_{nc}^+\;,
\end{equation}
where we assume that $(\alpha,\alpha)=2$ for all roots.
Suppose now that $\Lambda$ satisfies (\ref{first_condition}) and that $\alpha\in\Delta_{nc}^+$ is such that $\Lambda-\alpha$ is a weight of $V(\Lambda)$. Then the string $\{\Lambda-\alpha+j\beta\}$ extends to $-r\leq j\leq q$ with $r=q+ (\Lambda-\alpha,\beta)$. We are going to show that there exists $\beta\in\Delta^+_{nc}$ such that $(\Lambda-\alpha,\beta)=1$ so that
$r=q+1\geq 1$ and $e_{-\beta} e_{-\alpha} v_\Lambda\not = 0$. In particular this happens if there exists such $\beta$ such that $(\alpha,\beta)=0$ and $(\Lambda,\beta)=1$.

\bigskip
$\e_6$) Let $\R^8=\langle \epsilon_i, i=1,\ldots 8\rangle$. The real Cartan subalgebra $\t_\R=i\t$ can be described as the subspace of $\R^8$ generated by $\epsilon_i$, $i=1,\dots,5$ and $\epsilon=\epsilon_6+\epsilon_7+\epsilon_8$. As in \cite{Adams}, the positive roots $\Delta^+$ are
\begin{equation}\label{e6_roots}
\{\epsilon_{ij}=\epsilon_i- \epsilon_j\}_{1\leq i<j}^5
\coprod \{f_{ij}=(\epsilon_i+\epsilon_j)\}_{1\leq i<j}^5
\coprod \{\frac{1}{2}(\epsilon + \sum_{i=1}^5 s_i \epsilon_i), s_i=\pm, \Pi_i s_i=-1 \}\,.
\end{equation}
The simple roots are
\begin{equation}\label{e6_simple_roots}
\Pi=\{\alpha_i=\epsilon_{i,i+1}\}_{i=1}^4\coprod \{\alpha_5=f_{45},\alpha_6=\frac{1}{2}(\epsilon - \sum_{i=1}^5 \epsilon_i)\}
\end{equation}
and the non compact simple root $\phi=\alpha_6$. The non compact positive roots are
\begin{equation}\label{e6_noncompact_roots}
\Delta^+_{nc} = \{\phi + \frac{1}{2}\sum_{i=1}^5(1+s_i)\epsilon_i, \Pi_i s_i=-1\}\;.
\end{equation}
The weight $\Lambda = \sum_{i=1}^5 \Lambda_i \epsilon_i + \Lambda_0\epsilon$ is dominant if $(\Lambda,\alpha_i)=N_i$, $i=1,\ldots,6$, with $N_i\in\N$.
A straightforward computation gives
$$
\Lambda_1=N_1+N_2+N_3+\frac{1}{2}(N_4+N_5),\,\,\,\,\,\,
\Lambda_2 = N_2+N_3+\frac{1}{2}(N_4+N_5)\,,
$$
$$
\Lambda_3 = N_3+\frac{1}{2}(N_4+N_5)\,,\,\,
$$
$$
\Lambda_4 = \frac{1}{2}(N_4+N_5)\,,\,\,\, \Lambda_5 = \frac{1}{2}(N_5-N_4)\,,
$$
$$
\Lambda_0 = \frac{1}{3} N_1 +\frac{2}{3} N_2 + N_3 + \frac{1}{2}N_4 +\frac{5}{6} N_5 +\frac{2}{3} N_6\;.
$$

Condition (\ref{first_condition}) means
$$
(\Lambda,\phi) = N_6 \in\{0,1\},\,\, (\Lambda,\phi+\epsilon_i+\epsilon_j) = N_6+\Lambda_i+\Lambda_j\in\{0,1\},$$
$$(\Lambda,\phi+\epsilon_i+\epsilon_j+\epsilon_k+\epsilon_l)=N_6+\Lambda_i+\Lambda_j+\Lambda_k+\Lambda_l\in\{0,1\}\;,
$$
where $i,j,k,l=1,\ldots,5$. Let us consider first $N_6=1$. It is not difficult to verify that the only solution is $N_i=0$ for $i\not=6$ so that only $\Lambda =\frac{2}{3}\epsilon$ satisfies (\ref{first_condition}). Since $(\Lambda,\alpha)=1$ for each $\alpha\in\Delta_{nc}^+$, for every couple $\alpha,\beta$ of orthogonal non compact positive roots (which exists since the rank is $2$) $\Lambda-\alpha-\beta$ is a weight.

Let us consider now $N_6=0$. Since all couples $\Lambda_i+\Lambda_j$'s are non negative integers only one couple can be different from zero; and since $\Lambda_1\geq \Lambda_2 \geq \Lambda_3\geq \Lambda_4 \geq\Lambda_5$ it can be only $\Lambda_1+\Lambda_2$. This means $N_3=N_4=N_5=0$ and $\Lambda_1+\Lambda_2= N_1+2N_2 =1$ so that $N_1=1$ and $N_2=0$. The solution is then $\Lambda= \epsilon_1+\frac{1}{3}\epsilon$. Let us choose now $\alpha=\frac{1}{2}(\epsilon+\epsilon_1-\epsilon_2-\epsilon_3+\epsilon_4-\epsilon_5)$ and $\beta =\frac{1}{2}(\epsilon+\epsilon_1+\epsilon_2+\epsilon_3-\epsilon_4+\epsilon_5)$ so that $(\Lambda,\alpha)=(\Lambda,\beta)=1$ and $(\alpha,\beta)=0$ so that $\Lambda-\alpha-\beta$ is a weight.

\bigskip
$\e_7$) The real Cartan subalgebra $\t_\R=i\t$ is realized as the subspace of $\R^8$ generated by $\epsilon_i$, $1\leq i\leq 6$ and $\epsilon = \epsilon_7+\epsilon_8$. The simple roots can be chosen as
$$
\Pi = \{ \alpha_i=\epsilon_i-\epsilon_{i+1}, 1\leq i\leq 5;\, \alpha_6=\epsilon_5+\epsilon_6,\alpha_7 = \frac{1}{2}(\epsilon - \sum_{i=1}^6 \epsilon_i)\}\;
$$
and the non compact simple root $\phi=\alpha_1$ \cite{Adams}. The non compact positive roots $\Delta^+_{nc}$ are
$$
\{\epsilon_1-\epsilon_j,\epsilon_1+\epsilon_j; 2\leq j\leq 6\}\coprod \{\epsilon,\frac{1}{2}(\epsilon_1+\epsilon + \sum_{i=2}^6s_i \epsilon_i),s_i=\pm 1,\Pi_{i=2}^6s_i=1\}\;.
$$
The weight $\Lambda=\sum_{i=1}^6\Lambda_i \epsilon_i + \Lambda_0\epsilon$ is dominant if $(\Lambda,\alpha_i)=N_i\in\N$. A straightforward computation gives
$$
\Lambda_1=N_1+N_2+N_3+N_4+\frac{1}{2}(N_5+N_6),\,\,\,
\Lambda_2 = N_2+N_3+N_4+\frac{1}{2}(N_5+N_6)\,,$$
$$
\Lambda_3 = N_3+N_4+\frac{1}{2}(N_5+N_6)\,,\,\,
$$
$$
\Lambda_4 = N_4+\frac{1}{2}(N_5+N_6)\,,\,\,\, \Lambda_5 = \frac{1}{2}(N_5+N_6)\,, \Lambda_6 = \frac{1}{2}(N_6-N_5)
$$
$$
\Lambda_0 = N_0 +\frac{1}{2} N_1 +N_2+ \frac{3}{2} N_3 + 2N_4 + N_5+ \frac{3}{2}N_6\;.
$$
Condition (\ref{first_condition}) means
$$
(\Lambda,\epsilon_1+\epsilon_j) =\Lambda_1+\Lambda_j\in\{0,1\},\,\,
(\Lambda,\epsilon_1-\epsilon_j) = \Lambda_1-\Lambda_j\in\{0,1\},$$
$$
(\Lambda,\epsilon) = 2\Lambda_0\in\{0,1\},
$$
$$
(\Lambda,\frac{1}{2}(\epsilon_1+\epsilon +\sum_{i=1}^6s_i \epsilon_i)) = \frac{\Lambda_1}{2} + \Lambda_0 +\frac{1}{2} \sum_{i=2}^6s_i\Lambda_i\in\{0,1\},\,\,\;.
$$
Let us choose $N_1=(\Lambda,\epsilon_1-\epsilon_2)=0$. Then $(\Lambda,\epsilon_1+\epsilon_2)=\Lambda_1+\Lambda_2=2(N_2+N_3+N_4)+N_5+N_6\in\{0,1\}$ fixes $N_2=N_3=N_4=0$. Then $2\Lambda_0= 2N_0+2N_5 +3N_6\in\{0,1\}$ implies $N_0=N_5=N_6=0$ so that $\Lambda=0$.

Let us choose $N_1=1$. Then $\Lambda_1+\Lambda_2=1$ implies $N_2=\ldots=N_6=0$ and $2\Lambda_0=2N_0+1=1$ implies $N_0=0$. We then get $\Lambda = \epsilon_1 +\frac{1}{2}\epsilon$. It is clear that $(\Lambda,\alpha)=1$
for each non compact positive root $\alpha$, so that $\Lambda-\alpha-\beta$ is a weight for every couple $\alpha,\beta$ of orthogonal noncompact positive roots (the rank is $3$). \qed

\bigskip
\bigskip
\section{Invariant polynomials of $\k_\phi$}\label{invpol_section}
We are going to prove that, if we choose $\k_1=\k_\phi$ as defined in (\ref{main_subalgebra}), the condition (\ref{sufficient_condition}) is satisfied, so that, from Theorem \ref{thm_basic_forms}, the action of the Nijenhuis tensor on invariant $\k_\phi$-polynomials produces $\k_\phi$-basic forms.  We remark that $\k_\phi=\k_1$ in the Thimm chain defined in Section \ref{thimm_section} for all cases but $M_\phi=Gr(k,n)$ with $1<k<n-1$, {\it i.e.} Grassmanians that are not complex projective spaces.

The result is based on a local parametrization of the moment map around $\rho_\phi$. We refer to Section \ref{conventions} for notations and basic facts about compact hermitian symmetric spaces. Let us choose $P_\phi\subset\Delta^+_{nc}$ satisfying the following properties:
\begin{itemize}\label{maximal_noncompact_roots}
 \item[$i$)] for each $\alpha,\beta\in P_\phi$, $\alpha-\beta\not\in\Delta$; \item[$ii$)] $P_\phi$ is maximal with respect to ($i$).
\end{itemize}
Since the sum of two positive non compact roots is never a root, condition $i$) implies that roots in $P_\phi$ are all orthogonal and in particular linearly independent.

Let $\a_{P_\phi}\subset\k_\phi^\perp$ denote the space spanned by $\{i(e_\alpha+e_{-\alpha}),\alpha\in P_\phi\}$ and $\t_{P_\phi}=(P_\phi)^o\subset\t\subset\k_\phi$. We denote $\a_{P_\phi}'=J(\a_{P_\phi})$, that is the space spanned by $\{e_\alpha-e_{-\alpha},\alpha\in P_\phi\}$.

Let us consider the non compact real form $\g_0 = \k_\phi\dot{+}i\k_\phi^\perp$ of $\g$. We recall that a Cartan subalgebra $\h$ of $\g_0$ is maximally non compact if the non compact component $\h\cap i\k_\phi^\perp$ is a maximal abelian subalgebra of $\k_\phi^\perp$.

\begin{lemma}
$\h=\t_{P_\phi}\oplus i\a_{P_\phi}$ is a maximally non compact Cartan subalgebra of the non compact real form $\g_0$.
\end{lemma}
{\it Proof}. From the definitions it follows that $[\t_{P_\phi},\a_{P_\phi}]=0$ and $[\a_{P_\phi},\a_{P_\phi}]=0$. Since roots in $P_\phi$ are linearly independent, $\dim\t_{P_\phi}=\rk(\g)-\sharp P_\phi=\rk(\g)-\dim \a_{P_\phi}$. In particular $\h$ is a Cartan subalgebra. In order to prove that the non compact part $\a_{P_\phi}$ is maximal, let us suppose that there exists $\xi=\sum_{\alpha\in\Delta_{nc}^+} (\xi_\alpha e_\alpha+\xi_{-\alpha}e_{-\alpha})\in\k_\phi^\perp$ such that
$[\xi,a]=0$ for each $a\in\a_{P_\phi}$. Then for each $\beta\in P_\phi$ we have that
$$[\xi,e_\beta+e_{-\beta}]=(\xi_\beta-\xi_{-\beta}) [e_\beta,e_{-\beta}]+\sum_{\alpha\in \Delta_{nc}^+\setminus P_{\phi}}\xi_\alpha [e_\alpha,e_{-\beta}]+\xi_{-\alpha}[e_{-\alpha},e_\beta]=0\;\;,$$
where $[e_\beta,e_{-\beta}]\in i\t$.
First of all, we conclude that $\xi_\beta=\xi_{-\beta}$ for $\beta\in P_\phi$. Then, let us suppose that $\xi_\alpha\not=0$ for some $\alpha\in\Delta_{nc}^+\setminus P_\phi$; the above condition implies that $\alpha-\beta$ is not a root for all $\beta\in P_\phi$, but such root does not exist by maximality of $P_{\phi}$. We can then conclude that $\xi = \sum_{\beta\in P_\phi} \xi_\beta(e_\beta+e_{-\beta})\in \a_{P_\phi}$. \qed

\medskip
The compact Cartan is then decomposed as $\t=\t_{P_\phi}\oplus\t_{P_\phi}'$, where $\t_{P_\phi}'=\t_{P_\phi}^\perp\cap\t$. It is well known the $K_\phi$-orbits of any maximal abelian subalgebra of $\k_\phi^\perp$ cover all $\k_\phi^\perp$ (see Thm. 6.51 in \cite{Knapp}), {\it i.e.}
\begin{equation}\label{transitivity_non_compact_cartan}\k_\phi^\perp = \Ad_{K_\phi}(\a_{P_\phi}) \;.
\end{equation}
We need to prove the following properties.

\begin{lemma}\label{lemma3}
 \begin{itemize}
  \item[$i$)] $[\a_{P_\phi},\a_{P_\phi}']\subset \t_{P_\phi}'$;
  \item[$ii$)] $[\t,\a_{P_\phi}]=\a_{P_\phi}'$.
 \end{itemize}
\end{lemma}
{\it Proof}. Let us prove $i$). Let $\alpha,\beta\in P_\phi$ and compute
$$
[i(e_\alpha+e_{-\alpha}),e_\beta-e_{-\beta}] = -2i\delta_{\alpha,\beta} (e_\alpha,e_{-\alpha}) \tau_\alpha\in\t_{P_\phi}'
$$
where $(\tau_\alpha,H)=\alpha(H)$ for each $H\in\t$. Point $ii$) follows because $[H,i(e_\alpha+e_{-\alpha})]=i\alpha(H)(e_\alpha-e_{-\alpha})$ for each $H\in\t$. \qed

\medskip
The orthogonal decomposition $\k=\k_\phi\dot{+}\k_\phi^\perp$ assures that there exists an open neighborhood of $\rho_\phi\in M_\phi$ where the moment map $\mu\in C^\infty(M_\phi)\otimes\k$ can be written as
$$
\mu = \Ad_k \Ad_{e^\xi}(\rho_\phi)
$$
for $k\in K_\phi$ and $\xi\in\k_\phi^\perp$ so that $\mu\sim \Ad_{e^\xi}(\rho_\phi)$ where $\sim$ denotes up to $K_\phi$-adjoint action. By using (\ref{transitivity_non_compact_cartan}) and the fact that $\rho_\phi$ is $K_\phi$-invariant we can write
\begin{equation}\label{parametrization}
\mu \sim \Ad_{e^a}(\rho_\phi) \;.
\end{equation}
\begin{prop}
\label{general_corollary}
We have that
$$
[\mu_{\k_\phi},A_{\k_\phi}]=0
$$
where $A_{\k_\phi}$ is defined in (\ref{fundamental_object}). As a consequence, for each $p\in S(\k_\phi)^{\k_\phi}$, $d_Np$ is a basic $\k_\phi$-form.
\end{prop}
{\it Proof}.
Let us write
\begin{equation}\label{decomposition_non_compact_cartan}
\mu\sim \tilde{\rho}+\tilde{\xi}
\end{equation}
where $\tilde{\rho}\in\k_\phi$ and $\tilde{\xi}\in\k_\phi^\perp$.
From (\ref{parametrization}) we can write
$$
\tilde{\rho} = \sum_{k=0} \frac{1}{(2k)!} \ad_{a}^{2k}(\rho_\phi)\,,\;\;\;
\tilde{\xi} = \sum_{k=0} \frac{1}{(2k+1)!} \ad_{a}^{2k+1}(\rho_\phi)\;.
$$
It is clear that by using Lemma \ref{lemma3} $ii$) we see that $\ad_a(\rho_\phi)\in\a_{P_\phi}'$, from $i$) $\ad^2_a(\rho_\phi)\in\t_{P_\phi}'\subset\t$; we then show that $\ad^{2k}_a(\rho_\phi)\in\t$ and $\ad^{2k+1}_a(\rho_\phi)\in\a_{P_\phi}'$, so that
$\tilde\rho\in\t$ and $\tilde{\xi}\in\a_{P_\phi}'$. By recalling (\ref{fundamental_object}), we compute up to $K_\phi$-adjoint action that
$$
[\mu_{\k_\phi},A_{\k_\phi}] \sim \frac{1}{2} [\tilde{\rho},[J(\tilde{\xi}),\tilde{\xi}]] \in [\t,[J(\a_{P_\phi}'),\a_{P_\phi}']]=[\t,[\a_{P_\phi},\a_{P_\phi}']] = [\t,\t]=0\;,
$$
where we used Lemma \ref{lemma3} $i$). \qed

\medskip
We finally give the following explicit description of the terms appearing in (\ref{decomposition_non_compact_cartan}).

\begin{cor}
Let $P_\phi = \lbrace \alpha_j \rbrace_{j=1}^{\mathrm{rank}M_\phi}$, $\mathfrak{a}_{P_\phi}$ spanned by $X_{\alpha_j} =i(e_{\alpha_j} + e_{-\alpha_j}) $. 
Then, for any $a= \sum_{i=j}^{\mathrm{rank}M_\phi} a_j X_{\alpha_j}$ where $a_j\in \mathbb{R}$, we compute
\begin{equation}\label{resummation_1}
    \tilde{\rho}= \rho_\phi +  \sum_{j=1}^{\mathrm{rank}M_\phi} f_j \,ih_{\alpha_j},\qquad \tilde{\xi}= \sum_{j=1}^{\mathrm{rank}M_\phi} g_j JX_{\alpha_j},
\end{equation}
where
\begin{equation}\label{resummation_2}
    f_j= \frac{1}{2}( \cos (2a_j)-1),\qquad g_j= -\frac{1}{2}\sin (2a_j).
\end{equation}
Moreover,
\begin{equation}\label{A_resummation}
    [J\tilde{\xi}, \tilde{\xi}]= 2 \sum_{j=1}^{\mathrm{rank}M_\phi}\left[ f_j +  f^2_j \right] \,ih_{\alpha_j}.
\end{equation}
\end{cor}
{\it Proof}. First of all, let us remind that, for $\alpha_j,\alpha_l\in P_\phi$, $[ih_{\alpha_j},X_{\alpha_l}]=2 \delta_{jl}\,( JX_{\alpha_l})$ and $[X_{\alpha_j}, JX_{\alpha_l}]=2\delta_{jl}\, i h_{\alpha_l}$. Then, if we set
$$ \ad_a^{2k}(\rho_\phi)=\sum_j(f_j)_k \, ih_{\alpha_j},\qquad  \ad_a^{2k+1}(\rho_\phi)=\sum_j (g_j)_k JX_{\alpha_j},$$
we can find the following recursive relations
$$ (f_j)_k= -4 a_j^2 (f_j)_{k-1},\quad  (g_j)_k= -4 a_j^2 (g_j)_{k-1} $$
with $(f_j)_1=-2a_j^2$ and $(g_j)_0 = -a_j$. Then, we can write
$$ (f_j)_k = \frac{(-1)^k}{2} (2a_j)^{2k}, \quad (g_j)_k = \frac{(-1)^{k+1}}{2} (2a_j)^{2k+1}, $$
and so, putting $f_j = \sum_{k=1}^\infty \frac{1}{2k!}(f_j)_k $ and $g_j = \sum_{k=0}^\infty \frac{1}{(2k+1)!}(g_j)_k $, one obtains (\ref{resummation_1},\ref{resummation_2}). The last claim is obtained by a straightforward calculation and by considering that
$$ g^2_j = - f_j-f^2_j ,$$
formula obtained by (\ref{resummation_2}) using some basic goniometric relations. \qed

\medskip
Thanks to (\ref{resummation_1}) and (\ref{A_resummation}) $\k_\phi$- invariant polynomials can be expressed as polynomials in the variables $f_j$, $j=1\ldots,\rk(M_\phi) $. Let us denote with $p_n(f)=\sum_j c_j f_j^n$, where $c_j = i(h_{\alpha_j},\rho_\phi)=-2/(\alpha_j,\alpha_j)$. Remark that $(ih_{\alpha_j},ih_{\alpha_r})=2\delta_{jr} c_j$. Let us introduce the following $\k_\phi$-invariant polynomials
\begin{equation}\label{kphi_invariant}
I_{1,0} \equiv (\mu_{\k_\phi},\rho_\phi)= (\rho_\phi,\rho_\phi)+p_1\,,\;\; I_{2,0} \equiv (\mu_{\k_\phi},\mu_{\k_\phi})=(\rho_\phi,\rho_\phi)+2p_1+2p_2\,,
\end{equation}
$$
I_{0,1}=(\rho_\phi,A_{\k_\phi})= p_1+p_2 = \frac{1}{2}(I_{2,0}-(\rho_\phi,\rho_\phi))\,,\;\;\;I_{1,1} \equiv(\mu_{\k_\phi},A_{\k_\phi}) =  p_1+3p_2+2p_3\;.
$$
By applying (\ref{action_Nijenhuis}) we now easily compute
\begin{equation}
\label{Nijenhuis_kphi}
d_NI_{1,0} =d(I_{1,0}-\frac{1}{2} I_{2,0})\,,\;\;\; d_N I_{2,0} = d(I_{2,0}-\frac{2}{3} I_{1,0}- \frac{4}{3} I_{1,1})\;.
\end{equation}
Let us remark that
while $I_{0,1}\in \mu^*_{\k_\phi}(S(\k_\phi)^{\k_\phi})$, it is not obvious that the same is true for $I_{1,1}$ so that (\ref{action_Nijenhuis}) is not enough to compute $d_N I_{1,1}$. This property will be analyzed in Section \ref{exceptional} for the exceptional cases $EIII$ and $EVII$.

Let finally express the above formula on the basis given by the symmetric polynomials $p_j$. We get
\begin{equation}\label{Nijenhuis_kphi_symm}
d_Np_1  = -dp_2\, ,\;\;\;\; d_N p_2 = -\frac{4}{3} dp_3\;.
\end{equation}
It is tempting to state that
\begin{equation}
\label{Nijenhuis_eigenvalues_kphi}
d_Nf_i = -2 f_i df_i\;, i=1,\ldots \rk(M_\phi)
\end{equation}
that actually imply (\ref{Nijenhuis_kphi_symm}). Since in (\ref{Nijenhuis_eigenvalues_kphi}) there are two equations, if $\rk(M_\phi)=2$ then (\ref{Nijenhuis_eigenvalues_kphi}) and (\ref{Nijenhuis_kphi_symm}) are equivalent.  To see what happens when $\rk(M_\phi)>2$ it is useful to connect with the results of Theorem \ref{explicit_formula} and Corollary \ref{cor_indecomposable} for the classical $BCD$ cases.  Let us consider a $\phi\,$-minimal representation  $V(\Lambda)=V_+\oplus V_-$ of $\k$ (see Definition \ref{phi_minimal_rep}).  We need the following result.

\begin{lemma}
In a $\phi\,$-minimal representation, on $V_+$ we have that
$$
h_{\alpha_i}h_{\alpha_j} = \delta_{ij} h_{\alpha_i}
$$
for each $\alpha_i\in P_\phi$.
\end{lemma}

{\it Proof}. Let $\alpha$ be a positive non compact root and let us consider the $\sl(2,\C)$ subalgebra generated by $\{h_\alpha,e_\alpha,e_{-\alpha}\}$. Let $v_\lambda\in V_+$ be a common eigenvector of the Cartan subalgebra with
weight $\lambda$ such that $\lambda_\alpha=\lambda(h_\alpha)\not=0$. The eigenvector $v_{\lambda-\alpha}=e_{-\alpha}v_\lambda\in V_-$ is non vanishing (indeed $e_\alpha v_{\lambda-\alpha}=h_\alpha v_\lambda = \lambda_\alpha v_\lambda $). Let us consider
$W_\alpha= \langle v_\lambda,v_{\lambda-\alpha}\rangle$; since $e_{\alpha}v_\lambda=e_{-\alpha } v_{\lambda-\alpha}=0$, $W_\alpha$ is the two dimensional irreducible representation of $\sl(2,\C)$. As a consequence $\lambda_\alpha=1$. We can then conclude that the eigenvalues of $h_\alpha$ on $V_+$ can be $0,1$ (and so $h_\alpha$ is idempotent).

Let us denote with $\lambda_i=\lambda(h_{\alpha_i})$. For each $i\not = j$ we compute
$$
\lambda_i\lambda_j = ( v_\lambda, h_{\alpha_i}h_{\alpha_j}v_\lambda ) = ( v_\lambda, [e_{\alpha_i},e_{-\alpha_i}][e_{\alpha_j},e_{-\alpha_j}]v_\lambda)\;,
$$
where $(,)$ here denotes the scalar product that makes $V(\Lambda)$ unitary. 
Since $\alpha_i$ are non compact roots, in a $\phi\,$-minimal representation $e_{\alpha_i}e_{\alpha_j}=0$; since $\alpha_i\in P_\phi$ are mutually orthogonal then $[e_{\pm \alpha_i},e_{\mp \alpha_j}]=0$ for $i\not=j$. As a consequence $\lambda_i\lambda_j=0$. \qed

\medskip

We then conclude that the non zero eigenvalues of $\mu_{\k_\phi}$ in the representation $V_+$ are $\{\Lambda_\phi +i f_i\}_{i=1}^{\rk(M_\phi)}$. By applying the formula of Corollary \ref{cor_indecomposable} $ii$) we see that $-2f_i$ is a Nijenhuis eigenvalue, {\it i.e.} (\ref{Nijenhuis_eigenvalues_kphi}) holds.

\section{A formula for $S(\k_\phi)^{\k_\phi}$ on $EIII$ and $EVII$}\label{exceptional}
We are going to discuss here the exceptional cases. We saw in the previous section that minimal representations do not exist for $\e_6$ and $\e_7$ so that we cannot apply Corollary \ref{cor_indecomposable} to the Thimm chains discussed in Section \ref{thimm_section}. For $\k_\phi$, that is the first subalgebra of the chain, we can apply Proposition \ref{general_corollary}: we are going to compute the explicit form of the subcomplex of the Nijenhuis complex generated by $\k_\phi$-invariant polynomials.

\subsection{EIII}
Let us consider $M_\phi=E_6/SO(10)\times SO(2)$. We are going to give the explicit formulas for the action of the Nijenhuis tensor on $\so(10)\oplus \so(2)$-invariant polynomials.

By using the description of positive roots and non compact positive roots of $\e_6$ given in (\ref{e6_roots}) and (\ref{e6_noncompact_roots}) respectively, we get that
\begin{equation}\label{e6_P}
P_\phi = \left\lbrace \phi =\frac{1}{2}(\epsilon -\sum_{i=1}^5 \epsilon_i)  , \psi = \frac{1}{2}(\epsilon +\epsilon_1+\epsilon_2+\epsilon_3+\epsilon_4-\epsilon_5)  \right\rbrace \subset \Delta^+_{nc}
\end{equation}
are orthogonal noncompact roots. Since $\rk(M_\phi)=2$, this subset is maximal. Let $\a_{P_\phi}$ be the spaces spanned by $\{ X_\phi = i(e_\phi + e_{-\phi}),X_\psi = i(e_\psi + e_{-\psi})\}$, according to Section \ref{invpol_section}. Thus, formulas (\ref{resummation_1}) and (\ref{resummation_2}) are written as
\begin{equation}\begin{aligned}
    &\tilde{\rho} = \rho_\phi + f_\phi \, ih_{\phi} + f_\psi \, ih_{\psi},\\ &\tilde{\xi}= g_\phi JX_\phi +g_\psi JX_\psi. \end{aligned}
\end{equation}

Since $\rk(M_\phi)=2$ the formulas computed in (\ref{Nijenhuis_kphi}) completely describe the action of $d_N$ on $\so(10)\oplus \so(2)$-invariant polynomials. Indeed from (\ref{e6_P}) we have that $c_\phi=i(h_{\alpha_\phi},\rho_\phi)=c_\psi=i(h_{\alpha_\psi},\rho_\phi)=-1$ while $(\rho_\phi,\rho_\phi)=-\frac{4}{3}$. Then, $p_n (f)= - f_\phi^n - f_\psi^n$ and we have the following polynomial relation
$$
p_3= - \frac{3}{2}p_1p_2 - \frac{1}{2} p_1^3.
$$
From (\ref{Nijenhuis_kphi_symm}) we get that
$$
d_Np_1  = -dp_2\, ,\;\;\;\; d_N p_2 = \frac{4}{6} d\left( 3p_1p_2 + p_1^3\right)\;
$$
so that $p_1$ and $p_2$ generate the subcomplex $(\Omega^{\so(10)\oplus \so(2)}_{pol},d_N)$. We can conclude  that 
$$
d_N f_\phi = -2 f_\phi df_\phi\, ,\;\;\;\; d_N f_\psi = -2 f_\psi df_\psi \,,  
$$
{\it i.e.} $-2f_\phi$ and $-2f_\psi$ are Nijenhuis eigenvalues.

\subsection{EVII}
The case $M_\phi=E_7/E_6\times SO(2)$ presents a new feature. Indeed, since $\rk(M_\phi)=3$ we need three generators to describe $\mu_{\k_\phi}^*(S(\k_\phi))^{\k_\phi}$. Two of them can be $I_{1,0}=(\rho_\phi,\mu_{\k_\phi})$ and $I_{2,0}=(\mu_{\k_\phi},\mu_{\k_\phi})$, we have to look for the third one. The general formulas (\ref{Nijenhuis_kphi}) involve also $I_{1,1}=(\mu_{\k_\phi},A_{\k_\phi})$ which, by looking at the expression given in (\ref{Nijenhuis_kphi}),
as a cubic polynomials in the $f$'s, is clearly independent. The problem is now to understand if $I_{1,1}\in\mu_{\k_\phi}^*(S(\k_\phi)^{\k_\phi})$. In the previous case of $EIII$, since the rank was $2$, $I_{1,1}$, as a polynomial of (maximal) degree $3$ in two variables, was clearly generated by $I_{1,0}$ and $I_{2,0}$. The next independent generator in $S(\e_6)^{\e_6}$ appears in degree $5$; we can write it as $I_{5,0}=\Tr_{V_\Lambda}\mu_{\e_6}^5$ for some representation $V_\Lambda$ of $\e_6$; it is a polynomial of degree $5$ in the $f$-variables. By counting the $f$-degree, we easily conclude that $I_{1,1}$ cannot be polynomially generated by $\{I_{1,0},I_{2,0},I_{5,0}\}$ and $I_{5,0}$ is a polynomial function of $\{I_{1,0},I_{2,0},I_{1,1}\}$. By inverting this relation, $I_{1,1}$ can be expressed as a non polynomial function of the moment map, {\it i.e.} $I_{1,1}\in\mu_{\e_6}^*(C^\infty(\e_6^*)^{\e_6})$. We can then conclude that if we want to preserve the polynomial structure, we are forced to deal with the bigger ring that includes also polynomials in the variables $A_{\e_6}$. More importantly, the computation of $d_NI_{1,1}$ requires an extension of (\ref{action_Nijenhuis}) to this bigger ring of invariants.

\section{Conclusions}
In this paper we discussed two different approaches to the problem of diagonalizing the Nijenhuis tensor $N$ on the hermitian symmetric spaces $M_\phi=K/K_\phi$. The first one is based on the existence of a special representation of $\k$ that we call $\phi\,$-minimal: in this case we can define a Thimm chain of subalgebras $\k\supset \k_1\ldots\supset\k_i\ldots$ together with the set of generators $I_r$ of the ring of invariant polynomials $S(\k_i)^{\k_i}$ that make the diagonalization problem easily solved.

The second construction is based on a local parametrization around $\rho_\phi\in M_\phi$. This method is valid for all cases, including the exceptional ones, and produces the subset of the spectrum related to the subalgebra $\k_\phi$. It detects a different behavior of the EVII case: in fact the ring of invariant polynomials $S(\k_\phi)^{\k_\phi}$ is not preserved by the algebroid differential $d_N$ and the bigger ring of polynomials of $\mu_{\k_\phi}$ and $A_{\k_\phi}$ must be taken into consideration. Formula (\ref{action_Nijenhuis}) must be then generalized. Moreover, preliminary computations for the second subalgebra $\k_2=\so(8)\oplus\so(2)$ in the Thimm chain of $EIII$ show that the sufficient condition $[\mu_{\k_2},A_{\k_2}]=0$ in (\ref{sufficient_condition}) is not satisfied. These facts indicate that the cases $EIII$ and $EVII$ have an exceptional behavior. It is possible that it will be convenient to relate it to the fact that $S(\k)^\k\not= S(\a_{\k_\phi^\perp})^{\k_\phi}$ for the exceptional cases proved by Helgason in \cite{Helgason}. This problem will be addressed in a future publication.

\end{document}